\documentclass{elsart}               

\usepackage{graphicx}          
                               
\usepackage{comment}
\usepackage[normalem]{ulem}
\usepackage{multicol}
\usepackage{lipsum}

\newcounter{comments}
\usepackage{xcolor}

\definecolor{Red}{rgb}{0.8,0,0}
\definecolor{Green}{rgb}{0.2,0.6,0.2}
\definecolor{Blue}{rgb}{0,0,0.8}
\definecolor{Purple}{rgb}{.7,0,.7}

\usepackage{amsfonts}
\usepackage{amsmath}

\newtheorem{remark}{Remark}[section]
\newtheorem{theorem}{Theorem}[section]
\newtheorem{definition}{Definition}[section]
\newtheorem{corollary}{Corollary}[section]
\newtheorem{proposition}{Proposition}[section]
\newtheorem{example}{Example}[section]

\newcommand{\F}{\mathbb F_2}

\begin{document}

\begin{frontmatter}
\title{Modular Construction of Boolean Networks} 

\thanks[correspondingauthor]{Contributed equally}

\author[Florida]{Matthew Wheeler\thanksref{correspondingauthor}}\ead{mwheeler1@ufl.edu}
,    
\author[Iowa]{Claus Kadelka\thanksref{correspondingauthor}}\ead{ckadelka@iastate.edu},               
\author[Dayton]{Alan Veliz-Cuba}\ead{avelizcuba1@udayton.edu},               
\author[Kentucky]{David Murrugarra}\ead{murrugarra@uky.edu}, 
\author[Florida]{Reinhard Laubenbacher}\ead{reinhard.laubenbacher@medicine.ufl.edu}  

\address[Florida]{Department of Medicine, University of Florida, Gainesville, FL, United States}  
\address[Iowa]{Department of Mathematics, Iowa State University, Ames, IA, United States}             
\address[Dayton]{Department of Mathematics, University of Dayton, Dayton, OH, United States}        
\address[Kentucky]{Department of Mathematics, University of Kentucky,  Lexington, KY, United States}        
          

\begin{abstract}
Boolean networks have been used in a variety of settings, as models for general complex systems as well as models of specific systems in diverse fields, such as biology, engineering, and computer science. Traditionally, their properties as dynamical systems have been studied through simulation studies, due to a lack of mathematical structure. This paper uses a common mathematical technique to identify a class of  Boolean networks with a ``simple” structure and describes an algorithm to construct arbitrary extensions of a collection of simple Boolean networks. In this way, all Boolean networks can be obtained from a collection of simple Boolean networks as building blocks. The paper furthermore provides a formula for the number of extensions of given simple networks and, in some cases, provides a parametrization of those extensions. This has potential applications to the construction of networks with particular properties, for instance in synthetic biology, and can also be applied to develop efficient control algorithms for Boolean network models.
\end{abstract}


\end{frontmatter}

\section{Introduction}

Boolean networks have been used in a variety of settings, for theoretical studies as well as applications. A special class of Boolean networks, cellular automata, have served as models for complex systems and their emergent properties~\cite{wolfram1984cellular,schwab2020concepts}. General "edge of chaos" properties of biological systems have been studied using random Boolean networks~\cite{kauffman1993origins}. Their use as computational models of specific systems spans engineering, biology, computer science, and physics. This remarkable versatility arises from their structural simplicity as computational algorithms on the one hand, and their dynamic complexity as dynamical systems, on the other hand. The flip side of this dichotomy is that, when viewed as mathematical objects, there are few theoretical or algorithmic tools available to link their simple structure with their complex dynamics. The principal analysis tool is repeated simulation, that is, choosing many possible initializations for the computational algorithm and iterating it on each initialization until a limit cycle is reached. The result is overwhelmingly statistical properties obtaining by carrying out many simulation runs across large collections of Boolean networks, either random ones or chosen with specific properties. 

A Boolean network on $n$ variables can be thought of as a function from the set of binary strings of length $n$ to itself, described by $n$ Boolean coordinate functions that change the value of each of the variables, with input from some or all of the other variables. Iterative application of the function generates a time-discrete dynamical system on the set of binary strings of length $n$. It is not hard to see that \emph{any} such set function can be represented as a Boolean network, so that the class of Boolean networks is identical to the class of all set functions on binary strings. This dearth of mathematical structure explains why simulation is the analysis method of choice. However, in recent years, there have been attempts to bring mathematical tools to bear on this class of functions, with some success. An early example is~\cite{jarrah2007nested}, which presents a study of special classes of Boolean functions, using the fact that any Boolean function (hence, any function on binary strings) can be represented as a polynomial function over the field with two elements. This same approach led to the application of methods from computer algebra and algebraic geometry to reverse-engineer Boolean networks from observational data~\cite{veliz2012algebraic}. 

Finally, in~\cite{kadelka2023modularity}, we laid the foundation for a deeper mathematical theory of Boolean networks, following a common approach found across many fields of mathematics. The goal is to classify the axiomatically defined objects in that field, for instance finite groups or topological surfaces. The approach is philosophically similar to chemistry that has identified elementary chemical substances described by the periodic table, and has shown that every other chemical substance can be decomposed into these elements in specified ways. Several steps are involved: (i) identification of the ``simple" objects that cannot be further decomposed, (ii) a proof that any object can be decomposed into a collection of simple objects, (iii) a way to classify all simple objects, and (iv) an algorithm to build up objects from a collection of simple objects. In~\cite{kadelka2023modularity}, we addressed steps (i) and (ii). In this paper, we address step (iv), the assembly of simple objects into more complex ones. 

Section~\ref{sec:background} contains the basic background. 
In Section~\ref{sec:res_and_ext} we
introduce the framework needed to assemble networks. We first provide a definition of ``simple" Boolean networks, after which we describe an amalgamation construction that ``glues together" simple networks. 
Section~\ref{sec:enumeration} contains the main results. We show how one can explicitly parametrize and enumerate extensions for \textit{graphical} networks, which can be described completely by an edge-labeled graph. More generally, we enumerate the number of extensions for the class of \textit{nested canalizing} networks and arbitrary networks. The particular focus on networks governed by nested canalizing functions is due to their abundance in systems biology~\cite{kadelka2020meta}.
We conclude with a brief discussion of limitations of this study and open questions for future research. 

\section{Background}\label{sec:background}
In this section, we review some standard definitions and introduce the concept
of \emph{canalization}.
 Throughout the paper, let $\mathbb F_2$ be the binary field with elements 0 and 1. 

\subsection{Boolean functions}


The main class of Boolean functions that we consider in this paper will be nested-canalizing functions. 



\begin{definition}\label{def_essential}
A Boolean function $f:\F^n \to \F$ \emph{depends} on the variable $x_i$ if there exists an $\bf x \in \F^n$ such that 
$$f(\bf x)\neq f(\bf x + e_i),$$
where $e_i$ is the $i$th unit vector. In that case, we also say $f$ \emph{depends} on $x_i$. and $x_i$ is an \emph{essential} variable of $f$.
\end{definition}



\begin{definition}\label{def_canalizing}
A Boolean function $f:\F^n \to \F$ is \emph{canalizing} if there exists a variable $x_i$, $a,\,b\in\{0,\,1\}$ and a Boolean function $g(x_1,\ldots,x_{i-1},x_{i+1},\ldots,x_n)\not\equiv b$  such that
$$f(x_1,x_2,...,x_n)= \begin{cases}
b,& \ \text{if}\ x_i=a\\
g(x_1,x_2,...,x_{i-1},x_{i+1},...,x_n),& \ \text{if}\ x_i\neq a
\end{cases}$$
In that case, we say that $x_i$ \emph{canalizes} $f$ (to $b$) and call $a$ the \emph{canalizing input} (of $x_i$) and $b$ the \emph{canalized output}.
\end{definition}

\begin{definition}\label{def_nested_canalizing}
A Boolean function $f:\F^n \to \F$ is \emph{nested canalizing} with respect to the permutation $\sigma \in \mathcal{S}_n$, inputs $a_1,\ldots,a_n$ and outputs $b_1,\ldots,b_n$, if
\begin{equation*}f(x_{1},\ldots,x_{n})=
\left\{\begin{array}[c]{ll}
b_{1} & x_{\sigma(1)} = a_1,\\
b_{2} & x_{\sigma(1)} \neq a_1, x_{\sigma(2)} = a_2,\\
b_{3} & x_{\sigma(1)} \neq a_1, x_{\sigma(2)} \neq a_2, x_{\sigma(3)} = a_3,\\
\vdots  & \vdots\\
b_{n} & x_{\sigma(1)} \neq a_1,\ldots,x_{\sigma(n-1)}\neq a_{n-1}, x_{\sigma(n)} = a_n,\\
1 + b_{n} & x_{\sigma(1)} \neq a_1,\ldots,x_{\sigma(n-1)}\neq a_{n-1}, x_{\sigma(n)} \neq a_n.
\end{array}\right.\end{equation*}
The last line ensures that $f$ actually depends on all $n$ variables. From now on, we use the acronym NCF for nested canalizing function.
\end{definition}

We restate the following powerful stratification theorem for reference, which provides a unique polynomial form for any Boolean function.

\begin{theorem}[\cite{he2016stratification}]\label{thm:he}
Every Boolean function $f(x_1,\ldots,x_n)\not\equiv 0$ can be uniquely written as 
\begin{equation}\label{eq:matts_theorem}
    f(x_1,\ldots,x_n) = M_1(M_2(\cdots (M_{r-1}(M_rp_C + 1) + 1)\cdots)+ 1)+ q,
\end{equation}

where each $M_i = \displaystyle\prod_{j=1}^{k_i} (x_{i_j} + a_{i_j})$ is a nonconstant extended monomial, $p_C$ is the \emph{core polynomial} of $f$, and $k = \displaystyle\sum_{i=1}^r k_i$ is the canalizing depth. Each $x_i$ appears in exactly one of $\{M_1,\ldots,M_r,p_C\}$, and the only restrictions are the following ``exceptional cases'':
\begin{enumerate}
    \item If $p_C\equiv 1$ and $r\neq 1$, then $k_r\geq 2$;
    \item If $p_C\equiv 1$ and $r = 1$ and $k_1=1$, then $q=0$.
\end{enumerate}
When $f$ is not canalizing (\textit{i.e.}, when $k=0$), we simply have $p_C = f$.
\end{theorem}

From Equation~\ref{eq:matts_theorem}, we can directly derive an important summary statistic of NCFs.

\begin{definition}\label{def_layers} \emph{(\cite{kadelka2017influence})}
Given an NCF $f(x_1,\ldots,x_n)$ represented as in Equation~\ref{eq:matts_theorem}, we call the extended monomials $M_i$ the \emph{layers} of f and define the \emph{layer structure} of $f$ as the vector $(k_1,\ldots,k_r)$, which describes the number of variables in each layer. Note that $k_1+\cdots+k_r = n$ and that by exceptional case 1 in Theorem~\ref{thm:he}, $k_r\geq 2$. 
\end{definition}

Although the problem is \textbf{NP}-hard, there exist multiple algorithms for finding the layer structure of an NCF~\cite{dimitrova2022revealing}. 

\begin{example}
The Boolean functions 
\begin{align*}
    f(x_1,x_2,x_3,x_4)&=x_1\wedge (\neg x_2\vee (x_3 \wedge x_4)) = x_1\left[x_2\left[x_3x_4+1\right]+1\right],\\
    g(x_1,x_2,x_3,x_4)&=x_1\wedge (\neg x_2\vee x_3 \vee x_4) = x_1\left[x_2\left(x_3+1\right)\left(x_4+1\right)+1\right]
\end{align*}
are both nested canalizing. $f$ consists of three layers with layer structure $(1,1,2)$, while $g$ possesses only two layers and layer structure $(1,3)$.
\end{example}




\subsection{Boolean Networks}
A Boolean network refers to a collection of Boolean functions which map binary strings to binary strings.  In particular we have the following.
\begin{definition}\label{def_local_model}
  A \emph{Boolean network} is an $n$-tuple of \emph{coordinate functions} $F:=(f_1,\dots,f_n)$, where $f_i\colon \F^n\to \F$. Each function $f_i$ uniquely determines a map
\begin{equation*}\label{eqn:F_i}
  F_i\colon \F^n\longrightarrow \F^n,\qquad F_i\colon(x_1,\dots,x_n)\longmapsto (x_1,\dots,{f_i(x)},\dots,x_n),
\end{equation*}
where $x=(x_1,\dots,x_n)$. Every Boolean network defines a canonical map, where the functions are synchronously updated:
\[
F\colon \F^n\longrightarrow \F^n,\qquad F\colon(x_1,\dots,x_n)\longmapsto(f_1(x),\dots,f_n(x)).
\]
\end{definition}
In this paper, we only consider this canonical map, i.e., we only consider synchronously updated Boolean networks. Frequently, we use the terms \emph{linear} or \emph{nested canalizing network}, which simply means Boolean networks where all update rules are from the respective family of functions. Some variables of $F$ may have a constant update function (i.e., remain fixed forever). We refer to these variables as \emph{external parameters}.

To each network, we can associate a directed graph which encodes the structure of dependencies between each variable in the network.


\begin{definition}\label{def_wiring_diagram}
The \emph{wiring diagram} of a Boolean network $F=(f_1,\ldots,f_n): \F^n\rightarrow \F^n$ is the directed graph with vertices $x_1,\ldots,x_n$ and an edge from $x_i$ to $x_j$ if $f_j$ depends on $x_i$. 
\end{definition}

External parameters are exactly those nodes in a wiring diagram that do not possess any incoming edges. 

\begin{figure}
    \centering
    \includegraphics[scale=1.5]{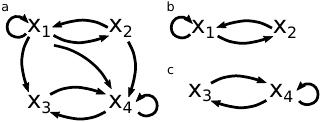}
    \caption{Boolean network decomposition. (a) Wiring diagram of a non-strongly connected Boolean network $F$. (b-c) Wiring diagram of $F$ restricted to (b) $\{x_1,x_2\}$ and (c) $\{x_3,x_4\}$. These strongly connected components are the wiring diagrams of the simple networks of $F$.
    }
    \label{fig:modules_eg}
\end{figure}

\begin{example}\label{eg:wd}
Figure~\ref{fig:modules_eg}a shows the wiring diagram of the Boolean network $F:\F^4\rightarrow \F^4$ given by \[F(x)=(x_2\wedge x_1, \neg x_1, x_1\vee \neg x_4, (x_1\wedge \neg x_2)\vee (x_3\wedge x_4)).\]
\end{example}

\section{Restrictions and extensions of Boolean functions and networks}\label{sec:res_and_ext}

In a previous paper, we derived two key results~\cite{kadelka2023modularity}. First, any Boolean network can be decomposed into a series of Boolean networks with strongly connected wiring diagrams, which we called modules in analogy to biology, with one-directional connections between the modules. Second, this structural decomposition implies a decomposition of the dynamics of the Boolean network. In this paper, we neglect dynamic considerations but dive deeper into structural questions. Specifically, the decomposition of a Boolean network into its modules can be considered as a kind of \emph{restriction}. On the contrary, two Boolean networks can be combined, with one upstream of the other, forming a kind of \emph{extension}. In this section, we formally define the restriction and extension of Boolean functions and Boolean networks, two inversely related concepts. 

\subsection{Restrictions and extensions of Boolean functions}\label{section:restriction}
In order to define restrictions of Boolean networks, we first define the concept of restrictions for Boolean functions.  

\begin{definition}
    Let $f:\F^n\to\F$ be a Boolean function and let $X=\{x_1,\ldots,x_n\}$ denote the variables of $f$. Let $Y\subset X$ be a subset of its variables.  A restriction of $f$ to $Y$, denoted by $f|_Y$, is a choice of values $x_i=a_i$ for each $x_i\not\in Y$. 
\end{definition}


\begin{example}
Consider the function
    \[f(x_1,x_2,x_3)=(x_1\wedge x_2)\vee(\neg x_1\wedge x_3)=\begin{cases}
        x_2 & \text{if }x_1=1,\\
        x_3 & \text{if }x_1=0.\\
    \end{cases}\]
In order to restrict $f$ to the set $Y=\{x_2,x_3\},$ we have two choices.  We may either set $x_1=1$ in which case $f|_Y=x_2$, or we can set $x_1=0$ yielding a different function $f|_Y=x_3.$ 
\end{example}

\begin{remark}\label{rem_unambiguous}
This example highlights that, in general, the restricted function depends on the specific choice of values $a_i$ and there exists no obvious default choice for the restriction. Certain classes of Boolean functions come, however, with such a natural, unambiguous choice of restriction. For nested canalizing functions, for example, there exists a natural choice of restriction: setting all variables not in $Y$ to their non-canalizing input values. This choice removes dependence on variables not in $Y$ without changing the way variables in $Y$ affect $f$.
\end{remark}



\begin{definition}\label{def_restr_ncf}
    Given a nested canalizing function $f:\F^n\to\F$ and a subset of its variables $Y\subset\{x_1,\ldots,x_n\},$ we define the restriction of $f$ to $Y$ to be the function such that for every $x_j\notin Y$ we set $x_j=a_j+1$ where $a_j$ is the canalizing input value for $x_j$.
    
    If $f$ depends only on a single variable, i.e., $f=x_j$ or $f=\neg x_j$, and $x_j \not\in Y$, then we set $f|_Y=1$. Note that this arbitrary choice of the constant agrees with the treatment of this case in exceptional case 1 in Theorem~\ref{thm:he}.
\end{definition}

\begin{example}
Given the nested canalizing function
    \[f(x_1,x_2,x_3)=\begin{cases}
                        b_{1} & x_1 = a_1,\\
b_{2} & x_1 \neq a_1, x_2 = a_2,\\
b_{3} & x_1 \neq a_1, x_2 \neq a_2, x_3 = a_3,\\
1 + b_{3} & x_1 \neq a_1, x_2 \neq a_2, x_3 \neq a_3,\\
                     \end{cases}
                     \]
and the subset $Y=\{x_1,x_3\},$ the restriction of $f$ to $Y$ is given by setting $x_2$ to its non-canalizing input,
    \[f|_Y(x_1,x_3)=f(x_1,a_2+1,x_3)=
    \begin{cases}
        b_{1} & x_1 = a_1,\\
        b_{3} & x_1 \neq a_1, x_3 = a_3,\\
        1 + b_{3} & x_1 \neq a_1, x_3 \neq a_3.
        \end{cases}\]
\end{example}

Having defined restrictions for a function, we can now introduce the notion of an extension of a Boolean function.
\begin{definition}
Let $X = \{x_1,\ldots,x_{n}\}, Y = \{y_1,\ldots,y_{q}\}$. An extension of a Boolean function $f:\F^{n}\to\F$ by $Y$ is any Boolean function $\tilde f:\F^{n+q}\to \F$ such that there exists a restriction $\tilde f|_X=f.$ That is, a restriction of an extension to the inputs of the original function recovers the original function.
\end{definition}

\begin{remark}
Note that some of the added variables of an extension may be non-essential. For example, $\tilde f(x,y) = x$ is an extension of $f(x) = x$.
\end{remark}

\begin{example}
Let $X = \{x_1,x_2\}, Y = \{y\}$ and $f = x_1 \oplus x_2 = x_1\ \text{XOR} \ x_2$ be a linear function (i.e., the exclusive OR function). Then, $g(x_1,x_2,y) = (x_1\ \text{XOR}\ x_2)\ \text{AND}\ y$ is an extension of $f$ because setting $y=1$ yields $g|_X = f$. Similarly, $g(x_1,x_2,y) = x_1\ \text{XOR}\ (x_2\ \text{AND}\ y)$ is an extension. On the other hand, $h(x_1,x_2,y) = (x_1\ \text{AND}\ y) \ \text{XOR}\ (x_2\ \text{OR}\ y)$ is not an extension. This is because setting $y=0$ yields $h|_X = x_2 \neq f$, and setting $y=1$ yields $h|_X = \text{NOT}\ x_1 \neq f$.
\end{example}

This example highlights that Boolean functions have many extensions. When considering Boolean functions in truth table format, the number of extensions of a Boolean function becomes evident.

\begin{theorem}\label{thm_n_ext_gen}
     Let $f:\F^{n}\to\F$ be a Boolean function. Then $f$ has 
     \[N_q(f) = \sum_{j=1}^{2^q} (-1)^{j+1} \binom{2^q}{j} 2^{2^{n+q}-j2^n} \] 
     extensions $\tilde f:\F^{n+q}\to \F$ by $q$ variables. Specifically, $f$ has $N_1(f) = 2^{2^{n}+1}-1$ extensions by a single variable.  
     
\end{theorem}

\begin{pf}
    We use the inclusion-exclusion principle to derive $N_q(f)$, the number of possible extensions of an $n$-variable Boolean function by $q$ variables. For each $\mathbf c = (c_1,\ldots,c_q)\in \F^q$, we define the set 
    $$A_{\mathbf c} = \big\{\tilde f: \F^{n+q}\to \F \big| \tilde f(x_1,\ldots,x_n,x_{n+1}=c_1,\ldots,x_{n+q}=c_q) = f(x_1,\ldots,x_n) \big\},$$
    which contains all extensions of $f$ that recover $f$ when the $q$ newly added variables are set to $\mathbf c$. In truth table format, any $\tilde f \in A_{\mathbf c}$ possesses $2^{n+q}$ rows, of which $2^n$ are required to match $f$. The other $2^{n+q} - 2^n$ rows can be freely chosen. Thus, $|A_{\mathbf c}| = 2^{2^{n+q}-2^n}$ for each $\mathbf c$.

    Some extensions $\tilde f$ recover $f$ for multiple choices of $\mathbf c \in \F^q$. For $m=1,\ldots,2^q$, the number of extensions that are in the intersection of $m$ different sets $A_{\mathbf c}$ can be computed with a similar argument as above. $m2^n$ rows of the truth table of $\tilde f$ are required to match $f$ for $m$ different choices of $\mathbf c \in \F^q$. The other $2^{n+q} - m2^n$ can be freely chosen. Thus, for any $C \subseteq \{1,\ldots,2^q\}$ with $|C| = m$, we have
    \[\Bigg|\bigcap_{\mathbf c \in C} A_{\mathbf c}\Bigg| = 2^{2^{n+q}-m2^n},\]
    which equals $1$ if $m=2^q$.
    Applying the inclusion-exclusion principle, we have
    \begin{align*}
    N_q(f) &= \Bigg|\bigcup_{\mathbf c \in \F^q} A_{\mathbf c}\Bigg|\\
    &= \binom{2^q}{1} |A_{\mathbf c}| - \binom{2^q}{2} |A_{\mathbf c_1} \cap A_{\mathbf c_2}| + \cdots + (-1)^{j+1}\binom{2^q}{j} |A_{\mathbf c_1} \cap \cdots \cap A_{\mathbf c_j}| + \cdots - 1\\
    &= \sum_{j=1}^{2^q} (-1)^{j+1} \binom{2^q}{j} 2^{2^{n+q}-j2^n}.
    \end{align*}

\end{pf}

\subsection{Restrictions and extensions of Boolean networks}
We next utilize the definitions for restriction and extension of functions to define similar concepts for Boolean networks. 
\begin{definition}
    Let $F=(f_1,\ldots,f_n):\F^n\to \F^n$ be a Boolean network. The network $$\tilde F = (\tilde f_1,\ldots,\tilde f_n,\tilde f_{n+1},\ldots,\tilde f_{n+m}):\F^{n+m}\to\F^{n+m}$$ is an extension of $F$ if for each $1\leq i\leq n$, $\tilde f_i$ is an extension of $f_i$.
\end{definition}
Then in a similar fashion as with functions, we can use the definition of a network extension to define a network restriction.
\begin{definition}
    Let $F:\F^{n+m}\to\F^{n+m}$ be a Boolean network. The network $G:\F^n\to\F^n$ is a restriction of $F$ if $F$ is an extension of $G$.
\end{definition}


\begin{example}
Consider the Boolean network 
\[F(x)=(x_2\wedge x_1, \neg x_1, x_1\vee \neg x_4, (x_1\wedge \neg x_2)\vee (x_3\wedge x_4))\] with wiring diagram in Figure~\ref{fig:modules_eg}a. The restriction of this network to $Y=\{x_1,x_2\}$ is the 2-variable network $F|_Y(x_1,x_2)=(x_2\wedge x_1, \neg x_1)$ with wiring diagram in Figure~\ref{fig:modules_eg}b, while the restriction of $F$ to $Y=\{x_3,x_4\}$ is the 2-variable network $F|_Y(x_3,x_4)=(\neg x_4, x_3\wedge x_4 )$ with wiring diagram in Figure~\ref{fig:modules_eg}c. Note that the wiring diagram of $F|_Y$ is always a subgraph of the wiring diagram of $F$, irrespective of the choice of $Y$.
\end{example}


\subsection{Simple Boolean networks}\label{subsec:simple_networks}
We now describe the elementary components in the structural decomposition theorem, which states that any Boolean network can be decomposed into a series of networks with strongly connected wiring diagrams, connected by a directed acyclic graph~\cite{kadelka2023modularity}. 


\begin{definition}\label{def_strongly_connected}
The wiring diagram of a Boolean network is \emph{strongly connected} if every pair of nodes is connected by a directed path. That is, for each pair of nodes $x_i,x_j$ in the wiring diagram  with $x_i\neq x_j$ there exists a directed path from $x_i$ to $x_j$ (and vice versa). In particular, a one-node wiring diagram is strongly connected by definition.
\end{definition}

\begin{example}
The wiring diagram in Figure~\ref{fig:modules_eg}a is not strongly connected but the wiring diagrams in Figure~\ref{fig:modules_eg}bc are strongly connected.
\end{example}

\begin{definition}\label{def_simple}
A Boolean network $F:\F^n\to \F^n$ is \emph{simple} if its wiring diagram is strongly connected. 
\end{definition}

\begin{remark}
In~\cite{kadelka2023modularity}, we defined a module of a network $F$ as a subnetwork of $F$ whose wiring diagram, excluding external inputs which it may or may not posses, is strongly connected. As such, every module of a network $F$ has an underlying simple network associated to it, which is just the restriction of the module determined by fixing the external inputs of the module to a specific value. 
\end{remark}


\begin{remark}
For an arbitrary Boolean network $F$, its wiring diagram will either be strongly connected or it will be made up of a collection of strongly connected components where connections between each component move in only one direction.  

Let $F$ be a Boolean network and let $W_1,\ldots,W_m$ be the strongly connected components of its wiring diagram, with $X_i$ denoting the set of variables in strongly connected component $W_i$. Then, $F$ can be decomposed into the \emph{simple networks}, defined by $F|_{X_1},\ldots,F|_{X_m}$.
\end{remark}

\begin{definition}\label{def:acyclic}
Let $W_1,\ldots, W_m$ be the strongly connected components of the wiring diagram of a Boolean network $F$. By setting $W_i \rightarrow W_j$ if there exists at least one edge from a vertex in $W_i$ to a vertex in $W_j$, we obtain a (directed) acyclic graph 
$$Q = \{(i,j) | W_i \rightarrow W_j \},$$ which describes the connections between the strongly connected components of $F$. 

\end{definition}


\begin{example}
For the Boolean network $F$ from Example~\ref{eg:wd}, the wiring diagram
has two strongly connected components $W_1$ and $W_2$ with variables $X_1=\{x_1,x_2\}$ and $X_2=\{x_3,x_4\}$ (Figure~\ref{fig:modules_eg}bc), connected according to the directed acyclic graph $Q = \{(1,2)\}$. The two simple networks of $F$ are given by the restriction of $F$ to $X_1$ and $X_2$, that is, $F|_{X_1}(x_1,x_2) = (x_2\wedge x_1, \neg x_1)$ and $F|_{X_2}(x_3,x_4) = (\neg x_4, x_3\wedge x_4)$. 
Note that the simple network $F|_{X_1}$, i.e., the restriction of $F$ to $X_1$, is simply the projection of $F$ onto the variables $X_1$ because $W_1$ does not receive feedback from the other component (i.e., because $(2,1)\not\in Q$).
\end{example}

\section{Parametrizing Network Extensions}\label{sec:enumeration}
The structural decomposition theorem shows how any Boolean network can be iteratively \emph{restricted} into a series of simple networks~\cite{kadelka2023modularity}. In this section, we consider the opposite: how and in how many ways can several (not necessarily simple) networks be combined into a larger network. For two networks, this can be thought of as determining the number of ways that the functions of one network can be \emph{extended} to include inputs coming from the other network. In general, this problem can be quite difficult. However, for certain families of functions with sufficient structure we can enumerate all such extensions.


\subsection{Graphical Boolean networks}
Various families of Boolean functions feature symmetries such that an $n$-node Boolean network governed by functions from such a family, including its dynamic update rules, is completely determined by a graph with $n$ nodes and labeled edges. We call networks governed by such a family of functions \emph{graphical}. The simplest graphical networks are those that are completely described by their wiring diagrams, which possesses two labels indicating presence or absence of a regulation. We call these networks \emph{2-graphical}. This includes Boolean networks governed entirely by linear, conjunctive (i.e., AND) or disjunctive (i.e., OR) functions. Boolean networks that can be completely described by a graph with 3 labels are called \emph{3-graphical} networks, and include those governed entirely by AND-NOT or OR-NOT functions. Here, the three labels indicate positive, negative or absent regulation.

\begin{definition}
Consider a class of Boolean networks governed by a family of Boolean functions such that the class is isomorphic to $\mathbb{F}_z^{n,n}$.
We call these networks \emph{z-graphical}. Further, we assume, without loss of generality, that the edge labels are $\mathbb F_z = \{0,1,\ldots,z-1\}$, with edge label $0$ ($\geq 1$) indicating absence (presence) of a regulation. 
\end{definition}

The set of all extensions of graphical networks can be enumerated and parameterized in a straightforward manner.




\begin{remark}\label{rem_partial_order}
Given an $n$-node $z$-graphical Boolean network $F$, the graph $\mathcal{G}$ that completely describes $F$ can be represented as an $n\times n$-square matrix with entries in $\mathbb F_z$. 
With the convention that edge label $0$ indicates absence of a regulation, $\mathcal{G}$ (just like the wiring diagram of $F$) can be represented by a lower block triangular matrix after a permutation of variables. The blocks on the diagonal represent the restrictions of $F$ to each simple network and the off-diagonal blocks represent uni-directional connections between the simple networks. The ordering of the blocks along the diagonal corresponds to the  acyclic graph $Q$ from Definition~\ref{def:acyclic}, which defines a partial ordering of the simple networks. Note that there exists only one block if $F$ is already simple.

Conversely, a collection of simple $z$-graphical networks can be combined to form a larger $z$-graphical network. The matrix representation of the graph that completely describes this larger network will consist of blocks along the diagonal representing the simple networks and blocks below the diagonal representing connections between each simple network.
\end{remark}

\begin{proposition}\label{prop:counting_graphical}
Let $X_1 = \F^{n_1}, X_2 = \F^{n_2}$. Given two $z$-graphical Boolean networks $F_1: X_1\to X_1,\,F_2: X_2\to X_2$, all possible graphs of an extension of $F_2$ by $F_1: X_1\times X_2\to X_1\times X_2$, can be parametrized by $\mathbb F_z^{n_2,n_1}$. 
\end{proposition}
\begin{pf}
Let $W_1, W_2$ be the square graphs in matrix form that completely describe $F_1$ and $F_2$. For a network $F$ to be an extension of $F_2$ by $F_1$, any $f_2\in F_2$ may depend on any $x_1\in X_1$ but the opposite is not admissible. The graph, in matrix form, of such an extension is thus given by
    \[W=\begin{bmatrix}W_1 & 0\\ \tilde P & W_2\end{bmatrix},\]
where the submatrix $\tilde P\in \mathbb F_z^{n_2,n_1}$  specifies the dependence of $X_2$ on $X_1$.
\end{pf}

\begin{corollary}\label{corr_number_extensions_graphical}
Let $X_1 = \F^{n_1}, X_2 = \F^{n_2}$. Given two z-graphical Boolean networks $F_1: X_1\to X_1,\,F_2: X_2\to X_2$, 
the number of different extensions of $F_2$ by $F_1: X_1\times X_2\to X_1\times X_2$ is $2^{n_1n_2}$.
\end{corollary}

\begin{example}
Consider the 2-graphical linear Boolean networks $F_1(x_1,x_2)=(x_1 + x_2,x_1)$ and $F_2(x_3,x_4)=(x_4,x_3+x_4)$, which are completely described by their wiring diagrams, shown in Figure~\ref{fig:modules_eg}bc. Figure~\ref{fig:modules_eg}a shows the wiring diagram of one of $2^4=16$ possible extensions of $F_2$ by $F_1$, 
$$F(x_1,x_2,x_3,x_4) = (x_1 + x_2,x_1,x_1+x_4,x_1+x_2+x_3+x_4).$$
In matrix form, we have
    \[F=\begin{bmatrix}F_1&0\\\tilde P&F_2\end{bmatrix},\; F_1=\begin{bmatrix}1&1\\1&0\end{bmatrix},\;F_2=\begin{bmatrix}0&1\\1&1\end{bmatrix},\;\tilde P=\begin{bmatrix}1&0\\1&1\end{bmatrix}.\]
The matrix $\tilde P$ is one out of a total of 16 Boolean $2\times2$ binary matrices, which determines the extension.
\end{example}

\begin{example}
Consider the 3-graphical AND-NOT networks 
\[F_1=(\neg x_2,x_1\wedge x_2)\quad\text{and}\quad F_2=(\neg x_3\wedge x_4,x_3),\]
which are completely described by the ternary $2\times 2$-matrices
\[F_1=\begin{bmatrix} 0&-1\\1&1\end{bmatrix},\quad F_2=\begin{bmatrix}-1 & 1\\1 & 0\end{bmatrix}.\]
Given the specific connection matrix $$\tilde P=\begin{bmatrix}0&0\\1&-1\end{bmatrix},$$ we can extend $F_2$ by $F_1$ to obtain the extended AND-NOT network
\[F= (\neg x_2,x_1\wedge x_2,\neg x_3\wedge x_4,x_1\wedge\neg x_2\wedge x_3),\]
which naturally decomposes into $F_1$ and $F_2$. Note that $\tilde P$ is one out of a total of $3^4 = 81$ ternary $2\times 2$ matrices, each of which determines a unique extension of $F_2$ by $F_1$.
\end{example}


We can generalize the previous results to any network, using the acyclic graph that encodes how the simple networks are connected, described in Remark~\ref{rem_partial_order}.



\begin{theorem}\label{thm:counting_graphical}
Let $X_1 = \F^{n_1}, \ldots, X_m = \F^{n_m}$ and let $F_1,\ldots,F_m$ with $F_i: X_i\to X_i$ be a collection of z-graphical Boolean networks. All possible graphs that completely describe a z-graphical network $F$, which admits a decomposition into simple networks $F_1,\ldots,F_m$, can be parametrized (i.e., are uniquely determined) by the choice of an acyclic graph $Q$ and an element in $$\prod_{\substack{(i,j)\in Q\\i\neq j}}\mathbb F_z^{n_j,n_i}\backslash \{0\}.$$ 
\end{theorem}
\begin{pf}
Let $F$ be a network with simple networks $F_1,\ldots,F_m$ and associated acyclic graph $Q$. (Note $F$ may already be simple). Further, let $W_i$ be the graphs, in matrix form, that completely describe $F_i, i=1,\ldots, m$. For each $(i,j) \in Q$, there exists a unique non-zero matrix $\tilde P_{ij} \in \mathbb F_z^{n_j,n_i}$, which describes the exact connections from $W_i$ to $W_j$. Note that the zero matrix would imply no connections from $W_i$ to $W_j$, contradicting $(i,j)\in Q$. 
\end{pf}

\begin{remark}
The set $\mathcal{Q}$ of all possible acyclic graphs corresponding to a Boolean network $F$ with $m\geq 1$ simple networks can be parametrized by the set of all lower triangular binary matrices with diagonal entries of 1. Thus,
$|\mathcal{Q}| = 2^{m(m-1)/2}$.
\end{remark}

\begin{corollary}
The number of all possible graphs completely describing a z-graphical Boolean network with simple networks $F_1,\ldots,F_m$ in that order (i.e., with any acyclic graph $Q$ satisfying $(i,j)\not\in Q$ whenever $i>j$) is given by
\[\sum_{Q\in\mathcal{Q}}\prod_{\substack{(i,j)\in Q\\i\neq j}}(z^{n_in_j}-1).\]
Note that the product over $Q=\emptyset$ is 1. Further, note that this number can also simply be expressed as $z^M$, where 
$$M = \sum_{i=1}^m\sum_{j=i+1}^m n_in_j$$
describes the total number of possible combinations of edge labels between all pairs of simple networks $F_i$ and $F_j$ with $i<j$.
\end{corollary}

\begin{corollary}
Let $X_1 = \F^{n_1}, \ldots, X_m = \F^{n_m}$ and let $F_1,\ldots,F_m$ with $F_i: X_i\to X_i$ be a collection of z-graphical Boolean networks. The number of different z-graphical networks with simple networks $F_1,\ldots,F_m$ with acyclic graph $Q$ is given by
$$\prod_{\substack{(i,j)\in Q\\i\neq j}}(z^{n_in_j}-1),$$
and any z-graphical network with simple networks $F_1,\ldots,F_m$ (with any  acyclic graph  $Q$ satisfying $(i,j)\not\in Q$ whenever $i>j$) is parameterized by an element of
\[\bigsqcup_{Q\in\mathcal{Q}}\prod_{\substack{(i,j)\in Q\\i\neq j}}\mathbb F_z^{n_j,n_i}\backslash \{0\}\cong \mathbb F_z^M\]
where 
$$M = \sum_{i=1}^m\sum_{j=i+1}^m n_in_j.$$
\end{corollary}

\subsection{Nested canalizing Boolean networks}
The results thus far enable us to parametrize and count all possible extensions for networks governed by families of Boolean functions that yield them completely described by a graph with edge labels. Biological Boolean network models are governed to a large part by NCFs~\cite{kadelka2020meta}. Boolean networks governed entirely by NCFs are not graphical since there are typically multiple possibilities how added variables can extend an NCF. In order to derive the number of possible extensions for networks governed by this important class of Boolean functions, we investigate what happens when we restrict an NCF one variable at a time. To simplify notation, let $\overline{\{x_i\}}$ denote the set of all variables of a function except for variable $x_i$.


\begin{example}\label{ex_ncf_count}
Consider the following NCF with unique polynomial representation
$$f=x_1(x_2+1)\left[x_3\left[(x_4+1)x_5+1\right]+1\right]+1$$
and layer structure $(2,1,2)$. Recall that, by Definition~\ref{def_restr_ncf}, restricting an NCF to $\overline{\{x_i\}}$ means setting $x_i$ to its non-canalizing value. We will now consider restrictions of $f$ to different $\overline{\{x_i\}}$ to highlight how the layer structure of $f$ determines the layer structure of the restriction.
\begin{enumerate}
    \item Restricting $f$ to $\overline{\{x_2\}}$ means setting $x_2=0$, i.e., removing the factor $(x_2+1)$ from the polynomial representation of $f$. This restriction has layer structure $(1,1,2)$ as the size of the first layer has been reduced by one.  
    \item Restricting $f$ to $\overline{\{x_3\}}$ means setting $x_3=1$ and yields
    $$f|_{\overline{\{x_3\}}} = x_1(x_2+1)(x_4+1)x_5+1,$$
    an NCF with a single layer. This reduction in the number of layers occurs because $x_3$ is the only variable in the second layer of $f$. Restricting to $\overline{\{x_3\}}$ implies an annihilation of this layer and results in a fusion of the adjacent layers (layer 1 and layer 3), as all variables in these two layers share the same canalized output value. 
    \item Restricting $f$ to $\overline{\{x_5\}}$ means setting $x_5=1$ and yields
    $$f|_{\overline{\{x_5\}}} = x_1(x_2+1)[x_3x_4+1]+1,$$
    which is an NCF of only two layers. This occurs because the last layer of an NCF always consists of two or more variables. Removal of $x_5$ from this layer leaves only variable $x_4$, which can now be combined with the previous layer. Note that the canalizing value of $x_4$ flips in this process.
\end{enumerate}
\end{example}
\begin{remark}\label{rem:addition_rules}
Given an NCF with layer structure $(k_1,\ldots,k_r)$, removing a variable in the $i$th layer (i.e., restricting the NCF to all but this variable) results in decreasing the size of the $i$th layer by one, $k_i\mapsto k_i-1.$ The previous example highlighted the three different cases that may occur when restricting NCFs. By reversing this thought process, we can categorize the number of ways that a new variable can be added to an NCF to obtain an extension.  
\begin{enumerate}
    \item (Initial layer) A new variable can always be added as new outermost layer. That is,
        \[(k_1,\ldots,k_r)\mapsto(1,k_1,\ldots,k_r)\]
    \item (Layer addition) A new variable can always be added to any layer. That is,
        \[(k_1,\ldots,k_r)\mapsto(k_1,\ldots,k_i+1,\ldots,k_r)\]
    \item (Splitting) If $k_i\geq 2,$ then the new variable can split the $i$th layer. That is,
        \[(k_1,\ldots,k_r)\mapsto(k_1,\ldots,k_i-\ell,1,\ell,\ldots,k_r)\]
        where $1\leq\ell\leq k_i-1$, with the exception that if $i=r$ and $\ell=1$,
        \[(k_1,\ldots,k_r)\mapsto(k_1,\ldots,k_r-1,2)\]
        because the last layer of an NCF always contains at least two variables.
\end{enumerate}
 Using these rules, we can construct extensions by adding new variables one at a time.
 \end{remark}
 
\begin{lem}\label{lem:ncf_unique_add}
Let $X=\F^{n_1}$. Consider an NCF $f:X\to\F$ and a set of variables $Y=\{y_1,\ldots,y_{n_2}\}$.  By fixing an ordering on $Y$, every extension $\tilde f: X\times Y\to\F$ of $f$ can be realized uniquely by the sequential addition of the variables of $Y$ using the rules described in Remark~\ref{rem:addition_rules}.
\end{lem}
\begin{pf}
The restriction of an NCF to $\overline{\{x_i\}}$ defines a unique function, as it means setting $x_i$ to its unique non-canalizing value. Let $\tilde f: X\times Y\to \F$ be an NCF extending $f$ over $Y$. By systematically restricting $\tilde f$ one variable at a time, we obtain a sequence of functions $\{f_i\}_{i=0}^{n_2}$
    \[f=f_0\to f_1\to \cdots\to f_i\to f_{i+1}\to \cdots f_{n_2}=\tilde f,\]
where each function $f_i$ is obtained by restricting $f_{i+1}$ to $\overline{\{y_i\}}.$ As every restriction of an NCF has an inverse extension (Definition~\ref{def_restr_ncf}), each $f_{i+1}$ can be obtained from $f_{i}$ by adding $y_{i}$ in the unique way that reverses the restriction.
\end{pf}

 \begin{figure}
    \centering
    \includegraphics{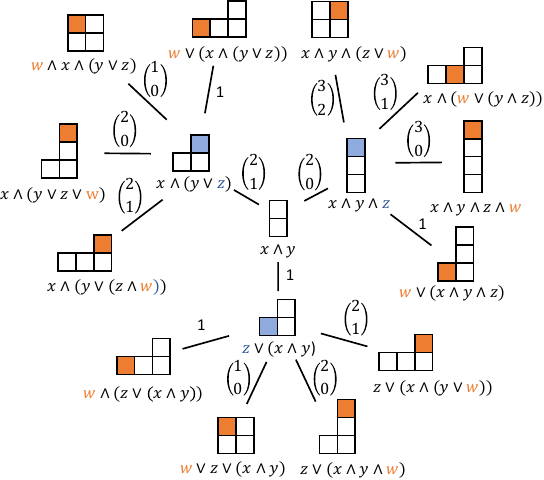}
    \caption{Graphical representation of the iterative process to find all possible nested canalizing extensions of a two-variable NCF. Each diagram represents a layer structure, each block represents a variable, each column represents a layer with the leftmost column representing the outer layer. The central diagram displays the only possible layer structure of a two-variable NCF. The colored blocks highlight the newly added variable; the first (second) added variable is shown in blue (orange) and within a layer the newly added variable is always shown on top. For each diagram, a representative nested canalizing extension is shown. Edge labels show the number of extensions sharing the given layer structure. Note that for each extended layer structure, the arbitrary choice of canalizing input value for the newly added variables supplies an additional factor of 2, which the edge labels do not include.}
    \label{fig:graphical_rep}
\end{figure} 

\begin{theorem}\label{thm:ncf_enumeration}
For an NCF $f$ with layer structure $(k_1,\ldots,k_r)$, the total number of unique extensions of $f$ by one variable is given by
    \[N_1(f) = 2+2\sum_{i=1}^r(2^{k_i}-1)\]
\end{theorem}
\begin{pf}
As outlined in Remark~\ref{rem:addition_rules}, there exist three distinct ways in which a variable $y$ can be added to an NCF to obtain a new, extended NCF: 1. as a new initial (i.e., outermost) layer, 2. as an addition to an existing layer, or 3. as a splitting of an existing layer. To obtain the total number of unique extensions of $f$ by one variable, we sum up all these possibilities, separately for each layer of $f$. Note that, irrespective of the way $y$ is added, there are always two choices for its canalizing value, which is why we multiply the total number by 2 in the end.

For a given layer of size $k_i$, there is exactly one way to add a single new variable $y$.  If $k_i=1$, the layer cannot be split so that  there is only one way to adjust this layer. If $k_i>1,$  the layer can be split
    \[(\ldots,k_i,\ldots)\mapsto (\ldots,k_i-\ell,1,\ell,\ldots),\]
and there are $\binom{k_i}{\ell}$ ways to choose $\ell$ of the $k_i$ variables to be placed in the layer following $y$.  A proper split requires $\ell,k_i-\ell\geq 1$ so that the number of possibilities to add $y$ to a layer (by adding to the layer or splitting it) is given by
\[M_i := 1+\sum_{j=1}^{k_i-1}\binom{k_i}{j}=2^{k_i}-1.\]
Note that this formula is also valid for the special last layer of an NCF, which always contains two or more variables: Splitting the final layer into a layer of size $k_r-1$ and a new last layer of size $\ell=1$ results in this last layer to be combined with the new layer consisting only of $y$, yielding a last layer of exact size $2$, which includes $y$ (see Case 3 in Example~\ref{ex_ncf_count} and Figure~\ref{fig:graphical_rep}).

Altogether, including the sole extension which adds $y$ as initial (i.e., outermost) later and accounting for the two choices of canalizing input value of $y$, the total number of unique extensions of $f$ by one variable is thus 
\[N_1(f) = 2\big(1+\sum_{i=1}^{r}M_i\big)=2+2\sum_{i=1}^{r}(2^{k_i}-1).\]
\end{pf}

\begin{remark}\label{rem_Nqf_NCF}
By combining Lemma~\ref{lem:ncf_unique_add} and Theorem~\ref{thm:ncf_enumeration}, one can enumerate all possible extensions of an NCF by a given number variables. Figure~\ref{fig:graphical_rep} illustrates the number of ways (up to choice of canalizing input value of the newly added variables) that two variables can be iteratively added to any NCF in two variables (which always consists of a single layer). Starting with a two-variable NCF $f$, there are
$$N_q(f) = 8, 92, 1328, 23184, 483840, 12050112$$
such extensions by $q=1,\ldots,6$ variables. This sequence did not have any match in the OEIS database~\cite{oeis}.

In general, all extensions of an NCF $f$ by $q$ variables, denoted $N_q(f)$, can be enumerated by an iterative application of Theorem~\ref{thm:ncf_enumeration}. 
\end{remark}


As for z-graphical Boolean networks (Corollary~\ref{corr_number_extensions_graphical}), we can compute in how many ways one nested canalizing Boolean network can be extended by another. 

\begin{proposition}\label{prop_ncf_network_ext}
Given two nested canalizing Boolean networks $F: \F^m\to \F^m, G = (g_1,\ldots,g_n): \F^n\to \F^n$, 
the number of different extensions of $G$ by $F: \F^{m+n}\to \F^{m+n}$ is 
\[\prod_{i=1}^n \sum_{Q\subseteq\{1,\ldots,m\}} N_{|Q|}(g_i) = \prod_{i=1}^n \sum_{q=0}^m \binom{m}{q} N_q(g_i),\]
where for all $g_i$, $N_q(g_i)$ is computed as described in Remark~\ref{rem_Nqf_NCF}, and $N_0(g_i) = 1$.
\end{proposition}

\begin{pf}
    Any combination of nodes in $F$ can be added as regulators to a node in $G$. Thus, for a given node in $G$ with update function $g_i$, the number of possible nested canalizing extensions of $g_i$ by the $m$ nodes in $F$ is given by 
    \[ \sum_{Q\subseteq\{1,\ldots,m\}} N_{|Q|}(g_i) = \sum_{q=0}^m \binom{m}{q} N_q(g_i),\]
    where (i) the sum iterates over all possible combinations of added regulators, and (ii) $N_q(g_i)$ can be derived as described in Remark~\ref{rem_Nqf_NCF}.
\end{pf}

\subsection{General Boolean networks}
Theorem~\ref{thm_n_ext_gen} provides a formula for the number of extensions of a general Boolean function, i.e., when not limiting ourselves to a specific family of functions. We can use the same formula as in Proposition~\ref{prop_ncf_network_ext} to count the number of Boolean network extensions, when any function is allowable.

\begin{proposition}\label{prop_gen_network_ext}
Given two Boolean networks $F: \F^m\to \F^m, G = (g_1,\ldots,g_n): \F^n\to \F^n$, 
the number of different extensions of $G$ by $F: \F^{m+n}\to \F^{m+n}$ is 
\[\prod_{i=1}^n \sum_{Q\subseteq\{1,\ldots,m\}} N_{|Q|}(g_i) = \prod_{i=1}^n \sum_{q=0}^m \binom{m}{q} N_q(g_i),\]
where for all $g_i$, $N_{q}(g_i)$ is computed as in Theorem~\ref{thm_n_ext_gen}, and $N_0(g_i) = 1$.
\end{proposition}

\section{Discussion}
Boolean networks play an important role as dynamical systems models for both theoretical and practical applications. In conjunction with our previous paper~\cite{kadelka2023modularity}, this article is laying the foundation for a mathematical structure theory for them that has implications for different applications and also points toward a future mathematical research program in this field. In~\cite{kadelka2023modularity}, we showed that every Boolean network can be decomposed into an amalgamation of simple networks. We showed that this structural decomposition induced a similar decomposition of their dynamics. As an application, we showed that control of Boolean networks can be done at the level of simple subnetworks, greatly simplifying the search for control interventions. 

In this paper, we considered the complementary problem of assembling Boolean networks from simple ones. 
This raises several technical challenges, some of which we addressed; others remain open. 
We enumerated all Boolean network extensions, as well as all extensions of specific classes of networks: z-graphical and nested canalizing ones.
Additionally for z-graphical networks, we obtained a straight-forward parametrization of all extensions. It remains to find a similar parametrization for all nested canalizing and general network extensions.
Other problems that remains are a categorical classification of all simple Boolean networks (similar to Hölder's classification of all finite simple groups), as well as the elucidation between their structure and their dynamics. As a related problem, it remains to be characterized when two extensions are dynamically equivalent. That is, all possible extensions of a collection of simple Boolean networks could be partitioned into dynamic equivalence classes, with a suitable definition of dynamic equivalence. 

Two other, potentially far-reaching, applications of the extension process can be envisioned. One is the construction of Boolean networks with specified properties by building them up from simple networks. One possible concrete setting for this might be synthetic biology. There, the genomes of complex organisms are assembled from simpler ``off-the-shelf" components. A related application is the construction of complex control systems from simple components with well-understood control mechanisms.

\section*{Acknowledgements}
Author Matthew Wheeler was supported by The American Association of Immunologists through an Intersect Fellowship for Computational Scientists and Immunologists. This work was further supported by the Simons foundation [grant numbers 712537 (to C.K.), 850896 (to D.M.), 516088 (to A.V.)]; the National Institute of Health [grant number 1 R01 HL169974-01 (to R.L.)]; and the Defense Advanced Research Projects Agency [grant number HR00112220038 (to R.L.)]. The authors also thank the Banff International Research Station for support through its Focused Research Group program during the week of May 29, 2022 (22frg001), which was of great help in framing initial ideas of this paper.

\bibliographystyle{plain}        
\bibliography{modularity_bib}           

\begin{thebibliography}{10}

\bibitem{dimitrova2022revealing}
Elena Dimitrova, Brandilyn Stigler, Claus Kadelka, and David Murrugarra.
\newblock Revealing the canalizing structure of {B}oolean functions: Algorithms
  and applications.
\newblock {\em Automatica}, 146:110630, 2022.

\bibitem{he2016stratification}
Qijun He and Matthew Macauley.
\newblock Stratification and enumeration of {B}oolean functions by canalizing
  depth.
\newblock {\em Physica D: Nonlinear Phenomena}, 314:1--8, 2016.

\bibitem{jarrah2007nested}
Abdul~Salam Jarrah, Blessilda Raposa, and Reinhard Laubenbacher.
\newblock Nested canalyzing, unate cascade, and polynomial functions.
\newblock {\em Physica D: Nonlinear Phenomena}, 233(2):167--174, 2007.

\bibitem{kadelka2020meta}
Claus Kadelka, Taras-Michael Butrie, Evan Hilton, Jack Kinseth, and Haris
  Serdarevic.
\newblock A meta-analysis of {B}oolean network models reveals design principles
  of gene regulatory networks.
\newblock {\em arXiv preprint arXiv:2009.01216}, 2020.

\bibitem{kadelka2017influence}
Claus Kadelka, Jack Kuipers, and Reinhard Laubenbacher.
\newblock The influence of canalization on the robustness of {B}oolean
  networks.
\newblock {\em Physica D: Nonlinear Phenomena}, 353:39--47, 2017.

\bibitem{kadelka2023modularity}
Claus Kadelka, Matthew Wheeler, Alan Veliz-Cuba, David Murrugarra, and Reinhard
  Laubenbacher.
\newblock Modularity of biological systems: a link between structure and
  function.
\newblock {\em Journal of the Royal Society Interface}, 20(207):20230505, 2023.

\bibitem{kauffman1993origins}
Stuart~A Kauffman.
\newblock {\em The origins of order: Self-organization and selection in
  evolution}.
\newblock Oxford University Press, USA, 1993.

\bibitem{oeis}
{OEIS Foundation Inc.}
\newblock The {O}n-{L}ine {E}ncyclopedia of {I}nteger {S}equences, 2024.
\newblock Published electronically at {http://oeis.org}.

\bibitem{schwab2020concepts}
Julian~D Schwab, Silke~D K{\"u}hlwein, Nensi Ikonomi, Michael K{\"u}hl, and
  Hans~A Kestler.
\newblock Concepts in {B}oolean network modeling: {W}hat do they all mean?
\newblock {\em Computational and structural biotechnology journal},
  18:571--582, 2020.

\bibitem{veliz2012algebraic}
Alan Veliz-Cuba.
\newblock An algebraic approach to reverse engineering finite dynamical systems
  arising from biology.
\newblock {\em SIAM Journal on Applied Dynamical Systems}, 11(1):31--48, 2012.

\bibitem{wolfram1984cellular}
Stephen Wolfram.
\newblock Cellular automata as models of complexity.
\newblock {\em Nature}, 311(5985):419--424, 1984.

\end{thebibliography}



\end{document}